# Bertrand Partner $D$-Curves in Minkowski 3-space $E_1^3$

**Mustafa Kazaz[a], H. Hüseyin Uğurlu[b], Mehmet Önder[a], Seda Oral[a]**

[a]*Celal Bayar University, Department of Mathematics, Faculty of Arts and Sciences, , Manisa, Turkey.*
E-mails:mustafa.kazaz@bayar.edu.tr, mehmet.onder@bayar.edu.tr
[b]*Gazi University, Gazi Faculty of Education, Department of Secondary Education Science and Mathematics Teaching, Mathematics Teaching Program, Ankara, Turkey.* E-mail: hugurlu@gazi.edu.tr

**Abstract**
In this paper, we consider the idea of Bertrand curves for curves lying on surfaces in Minkowski 3-space $E_1^3$. By considering the Darboux frame, we define these curves as Bertrand $D$-curves and give the characterizations for those curves. We also find the relations between the geodesic curvatures, the normal curvatures and the geodesic torsions of these associated curves. Furthermore, we show that in Minkowski 3-space $E_1^3$, the definition and the characterizations of Bertrand $D$-curves include those of Bertrand curves in some special cases.

## 1. Introduction

In the theory of space curves in differential geometry, the associated curves, the curves for which at the corresponding points of them one of the Frenet vectors of a curve coincides with the one of the Frenet vectors of the other curve have an important role for the characterizations of space curves. The well-known examples of such curves are Bertrand curves. These special curves are very interesting and characterized as a kind of corresponding relation between two curves such that the curves have the common principal normal i.e., the Bertrand curve is a curve which shares the normal line with another curve. These curves have an important role in the theory of curves. Hereby, from the past to today, a lot of mathematicians have studied on Bertrand curves in different areas [1,2,3,4,14,15,19]. Also these curves have an important role in the theory of ruled surface and in the characterizations of some other special curves. In [7], Izumiya and Takeuchi have studied cylindrical helices and Bertrand curves from the view point as curves on ruled surfaces. They have shown that cylindrical helices can be constructed from plane curves and Bertrand curves can be constructed from spherical curves. Also, they have studied generic properties of cylindrical helices and Bertrand curves as applications of singularity theory for plane curves and spherical curves[8]. In [3], Gluck has investigated the Bertrand curves in $n$-dimensional Euclidean space $E^n$. The corresponding characterizations of the Bertrand curves in $n$-dimensional Lorentzian space $E_1^n$ have been given by Tosun and Ozgür[15].

Furthermore, by considering the Frenet frame of the ruled surfaces, Ravani and Ku extended the notion of Bertrand curve to the ruled surfaces and named as Bertrand offsets[13]. The corresponding characterizations of the Bertrand offsets of timelike ruled surface in Minkowski 3-space $E_1^3$ were given by Kurnaz [10].

The differential geometry of the curves fully lying on a surface in Minkowski 3-space $E_1^3$ has been given by Ugurlu, Kocayigit and Topal[9,16,17,18]. They have given the Darboux frame of the curves according to the Lorentzian characters of surfaces and the curves.

In this paper, we consider the notion of the Bertrand curve for the curves lying on the surfaces in Minkowski 3-space $E_1^3$. We call these new associated curves as Bertrand $D$-curves and by using the Darboux frame of the curves we give the definition and the characterizations of these special curves.

## 2. Preliminaries

The Minkowski 3-space $E_1^3$ is the real vector space $IR^3$ provided with the standart flat metric given by
$$\langle,\rangle = -dx_1^2 + dx_2^2 + dx_3^2$$
where $(x_1, x_2, x_3)$ is a rectangular coordinate system of $E_1^3$. An arbitrary vector $\vec{v} = (v_1, v_2, v_3)$ in $E_1^3$ can have one of three Lorentzian causal characters; it can be spacelike if $\langle \vec{v}, \vec{v} \rangle > 0$ or $\vec{v} = 0$, timelike if $\langle \vec{v}, \vec{v} \rangle < 0$ and null (lightlike) if $\langle \vec{v}, \vec{v} \rangle = 0$ and $\vec{v} \neq 0$. Similarly, an arbitrary curve $\vec{\alpha} = \vec{\alpha}(s)$ can locally be spacelike, timelike or null (lightlike), if all of its velocity vectors $\alpha'(s)$ are respectively spacelike, timelike or null (lightlike). We say that a timelike vector is future pointing or past pointing if the first compound of the vector is positive or negative, respectively. For any vectors $\vec{x} = (x_1, x_2, x_3)$ and $\vec{y} = (y_1, y_2, y_3)$ in $E_1^3$, in the meaning of Lorentz vector product of $\vec{x}$ and $\vec{y}$ is defined by

$$\vec{x} \times \vec{y} = \begin{vmatrix} e_1 & -e_2 & -e_3 \\ x_1 & x_2 & x_3 \\ y_1 & y_2 & y_3 \end{vmatrix} = (x_2 y_3 - x_3 y_2, x_1 y_3 - x_3 y_1, x_2 y_1 - x_1 y_2)$$

where
$$\delta_{ij} = \begin{cases} 1 & i = j, \\ 0 & i \neq j, \end{cases} \quad e_i = (\delta_{i1}, \delta_{i2}, \delta_{i3}) \text{ and } e_1 \times e_2 = -e_3, \; e_2 \times e_3 = e_1, \; e_3 \times e_1 = -e_2.$$

Denote by $\{\vec{T}, \vec{N}, \vec{B}\}$ the moving Frenet frame along the curve $\alpha(s)$ in the Minkowski space $E_1^3$. For an arbitrary spacelike curve $\alpha(s)$ in the space $E_1^3$, the following Frenet formulae are given,

$$\begin{bmatrix} \vec{T}' \\ \vec{N}' \\ \vec{B}' \end{bmatrix} = \begin{bmatrix} 0 & k_1 & 0 \\ -\varepsilon k_1 & 0 & k_2 \\ 0 & k_2 & 0 \end{bmatrix} \begin{bmatrix} \vec{T} \\ \vec{N} \\ \vec{B} \end{bmatrix},$$

where $\langle \vec{T}, \vec{T} \rangle = 1$, $\langle \vec{N}, \vec{N} \rangle = \varepsilon = \pm 1$, $\langle \vec{B}, \vec{B} \rangle = -\varepsilon$, $\langle \vec{T}, \vec{N} \rangle = \langle \vec{T}, \vec{B} \rangle = \langle \vec{N}, \vec{B} \rangle = 0$ and $k_1$ and $k_2$ are curvature and torsion of the spacelike curve $\alpha(s)$ respectively. Here, $\varepsilon$ determines the kind of spacelike curve $\alpha(s)$. If $\varepsilon = 1$, then $\alpha(s)$ is a spacelike curve with spacelike first principal normal $\vec{N}$ and timelike binormal $\vec{B}$. If $\varepsilon = -1$, then $\alpha(s)$ is a spacelike curve with timelike principal normal $\vec{N}$ and spacelike binormal $\vec{B}$. Furthermore, for a timelike curve $\alpha(s)$ in the space $E_1^3$, the following Frenet formulae are given in as follows,

$$\begin{bmatrix} \vec{T}' \\ \vec{N}' \\ \vec{B}' \end{bmatrix} = \begin{bmatrix} 0 & k_1 & 0 \\ k_1 & 0 & k_2 \\ 0 & -k_2 & 0 \end{bmatrix} \begin{bmatrix} \vec{T} \\ \vec{N} \\ \vec{B} \end{bmatrix}.$$

where $\langle \vec{T}, \vec{T} \rangle = -1$, $\langle \vec{N}, \vec{N} \rangle = \langle \vec{B}, \vec{B} \rangle = 1$, $\langle \vec{T}, \vec{N} \rangle = \langle \vec{T}, \vec{B} \rangle = \langle \vec{N}, \vec{B} \rangle = 0$ and $k_1$ and $k_2$ are curvature and torsion of the timelike curve $\alpha(s)$ respectively[16,17].

**Definition 2.1. i) Hyperbolic angle:** Let $\vec{x}$ and $\vec{y}$ be future pointing (or past pointing) timelike vectors in $IR_1^3$. Then there is a unique real number $\theta \geq 0$ such that $<\vec{x}, \vec{y}> = -|\vec{x}||\vec{y}|\cosh\theta$. This number is called the *hyperbolic angle* between the vectors $\vec{x}$ and $\vec{y}$.

***ii) Central angle:*** Let $\vec{x}$ and $\vec{y}$ be spacelike vectors in $IR_1^3$ that span a timelike vector subspace. Then there is a unique real number $\theta \geq 0$ such that $<\vec{x}, \vec{y}> = |\vec{x}||\vec{y}|\cosh\theta$. This number is called the *central angle* between the vectors $\vec{x}$ and $\vec{y}$.

***iii) Spacelike angle:*** Let $\vec{x}$ and $\vec{y}$ be spacelike vectors in $IR_1^3$ that span a spacelike vector subspace. Then there is a unique real number $\theta \geq 0$ such that $<\vec{x}, \vec{y}> = |\vec{x}||\vec{y}|\cos\theta$. This number is called the *spacelike angle* between the vectors $\vec{x}$ and $\vec{y}$.

***iv) Lorentzian timelike angle:*** Let $\vec{x}$ be a spacelike vector and $\vec{y}$ be a timelike vector in $IR_1^3$. Then there is a unique real number $\theta \geq 0$ such that $<\vec{x}, \vec{y}> = |\vec{x}||\vec{y}|\sinh\theta$. This number is called the *Lorentzian timelike angle* between the vectors $\vec{x}$ and $\vec{y}$ [12].

**Definition 2.2.** A surface in the Minkowski 3-space $IR_1^3$ is called a timelike surface if the induced metric on the surface is a Lorentz metric and is called a spacelike surface if the induced metric on the surface is a positive definite Riemannian metric, i.e., the normal vector on the spacelike (timelike) surface is a timelike (spacelike) vector, [12].

**Lemma 2.1.** *In the Minkowski 3-space $IR_1^3$, the following properties are satisfied:*

*(i) Two timelike vectors are never orthogonal.*
*(ii) Two null vectors are orthogonal if and only if they are linearly dependent.*
*(iii) A timelike vector is never orthogonal to a null (lightlike) vector* [12].

## 3. Darboux Frame of a Curve Lying on a Surface in Minkowski 3-space $E_1^3$

Let $S$ be an oriented surface in three-dimensional Minkowski space $E_1^3$ and let consider a non-null curve $x(s)$ lying on $S$ fully. Since the curve $x(s)$ is also in space, there exists Frenet frame $\{\vec{T}, \vec{N}, \vec{B}\}$ at each points of the curve where $\vec{T}$ is unit tangent vector, $\vec{N}$ is principal normal vector and $\vec{B}$ is binormal vector, respectively.

Since the curve $x(s)$ lies on the surface $S$ there exists another frame of the curve $x(s)$ which is called Darboux frame and denoted by $\{\vec{T}, \vec{g}, \vec{n}\}$. In this frame $\vec{T}$ is the unit tangent of the curve, $\vec{n}$ is the unit normal of the surface $S$ and $\vec{g}$ is a unit vector given by $\vec{g} = \vec{n} \times \vec{T}$. Since the unit tangent $\vec{T}$ is common in both Frenet frame and Darboux frame, the vectors $\vec{N}$, $\vec{B}$, $\vec{g}$ and $\vec{n}$ lie on the same plane. Then, if the surface $S$ is an oriented timelike surface, the relations between these frames can be given as follows

If the curve $x(s)$ is timelike.

$$\begin{bmatrix} \vec{T} \\ \vec{g} \\ \vec{n} \end{bmatrix} = \begin{bmatrix} 1 & 0 & 0 \\ 0 & \cos\varphi & \sin\varphi \\ 0 & -\sin\varphi & \cos\varphi \end{bmatrix} \begin{bmatrix} \vec{T} \\ \vec{N} \\ \vec{B} \end{bmatrix},$$

If the curve $x(s)$ is spacelike.

$$\begin{bmatrix} \vec{T} \\ \vec{g} \\ \vec{n} \end{bmatrix} = \begin{bmatrix} 1 & 0 & 0 \\ 0 & \cosh\varphi & \sinh\varphi \\ 0 & \sinh\varphi & \cosh\varphi \end{bmatrix} \begin{bmatrix} \vec{T} \\ \vec{N} \\ \vec{B} \end{bmatrix}.$$

If the surface $S$ is an oriented spacelike surface, then the curve $x(s)$ lying on $S$ is a spacelike curve. So, the relations between the frames can be given as follows

$$\begin{bmatrix} \vec{T} \\ \vec{g} \\ \vec{n} \end{bmatrix} = \begin{bmatrix} 1 & 0 & 0 \\ 0 & \cosh\varphi & \sinh\varphi \\ 0 & \sinh\varphi & \cosh\varphi \end{bmatrix} \begin{bmatrix} \vec{T} \\ \vec{N} \\ \vec{B} \end{bmatrix}.$$

In all cases, $\varphi$ is the angle between the vectors $\vec{g}$ and $\vec{N}$.

According to the Lorentzian causal characters of the surface $S$ and the curve $x(s)$ lying on $S$, the derivative formulae of the Darboux frame can be changed as follows:

**i)** If the surface $S$ is a timelike surface, then the curve $x(s)$ lying on $S$ can be a spacelike or a timelike curve. Thus, the derivative formulae of the Darboux frame of $x(s)$ is given by

$$\begin{bmatrix} \dot{\vec{T}} \\ \dot{\vec{g}} \\ \dot{\vec{n}} \end{bmatrix} = \begin{bmatrix} 0 & k_g & -\varepsilon k_n \\ k_g & 0 & \varepsilon \tau_g \\ k_n & \tau_g & 0 \end{bmatrix} \begin{bmatrix} \vec{T} \\ \vec{g} \\ \vec{n} \end{bmatrix}, \quad \langle \vec{T},\vec{T}\rangle = \varepsilon = \pm 1, \ \langle \vec{g},\vec{g}\rangle = -\varepsilon, \ \langle \vec{n},\vec{n}\rangle = 1. \tag{$1_a$}$$

**ii)** If the surface $S$ is a spacelike surface, then the curve $x(s)$ lying on $S$ is a spacelike curve. Thus, the derivative formulae of the Darboux frame of $x(s)$ is given by

$$\begin{bmatrix} \dot{\vec{T}} \\ \dot{\vec{g}} \\ \dot{\vec{n}} \end{bmatrix} = \begin{bmatrix} 0 & k_g & k_n \\ -k_g & 0 & \tau_g \\ k_n & \tau_g & 0 \end{bmatrix} \begin{bmatrix} \vec{T} \\ \vec{g} \\ \vec{n} \end{bmatrix}, \quad \langle \vec{T},\vec{T}\rangle = 1, \ \langle \vec{g},\vec{g}\rangle = 1, \ \langle \vec{n},\vec{n}\rangle = -1. \tag{$1_b$}$$

In these formulae $k_g$, $k_n$ and $\tau_g$ are called the geodesic curvature, the normal curvature and the geodesic torsion, respectively. Here and in the following, we use "dot" to denote the derivative with respect to the arc length parameter of a curve.

The relations between geodesic curvature, normal curvature, geodesic torsion and $\kappa$, $\tau$ are given as follows

$$k_g = \kappa \cos\varphi, \ k_n = \kappa \sin\varphi, \ \tau_g = \tau + \frac{d\varphi}{ds}, \text{ if both } S \text{ and } x(s) \text{ are timelike or spacelike}, \tag{$2_a$}$$

$$k_g = \kappa \cosh\varphi, \ k_n = \kappa \sinh\varphi, \ \tau_g = \tau + \frac{d\varphi}{ds}, \text{ if } S \text{ is timelike and } x(s) \text{ is spacelike}. \tag{$2_b$}$$

(See [9,18,19]). Furthermore, the geodesic curvature $k_g$ and geodesic torsion $\tau_g$ of the curve $x(s)$ can be calculated as follows

$$k_g = \left\langle \frac{dx}{ds}, \frac{d^2x}{ds^2} \times \vec{n} \right\rangle, \ \tau_g = \left\langle \frac{dx}{ds}, \vec{n} \times \frac{d\vec{n}}{ds} \right\rangle \tag{3}$$

In the differential geometry of surfaces, for a curve $x(s)$ lying on a surface $S$ the followings are well-known

    **i)** $x(s)$ is a geodesic curve $\Leftrightarrow k_g = 0$,

    **ii)** $x(s)$ is an asymptotic line $\Leftrightarrow k_n = 0$,

    **iii)** $x(s)$ is a principal line $\Leftrightarrow \tau_g = 0$ [11].

## 4. Bertrand $D$-Curves in Minkowski 3-space $E_1^3$

In this section, by considering the Darboux frame, we define Bertrand $D$-curves and give the characterizations of these curves in Minkowski 3-space $E_1^3$.

**Definition 4.1.** Let $S$ and $S_1$ be oriented surfaces in Minkowski 3-space $E_1^3$ and let consider the arc-length parameter curves $x(s)$ and $x_1(s_1)$ lying fully on $S$ and $S_1$, respectively. Denote the Darboux frames of $x(s)$ and $x_1(s_1)$ by $\{\vec{T},\vec{g},\vec{n}\}$ and $\{\vec{T}_1,\vec{g}_1,\vec{n}_1\}$, respectively. If there exists a corresponding relationship between the curves $x$ and $x_1$ such that, at the corresponding points of the curves, the Darboux frame element $\vec{g}$ of $x$ coincides with the Darboux frame element $\vec{g}_1$ of $x_1$, then $x$ is called a Bertrand $D$-curve, and $x_1$ is a Bertrand partner $D$-curve of $x$. Then, the pair $\{x, x_1\}$ is said to be a Bertrand $D$-pair. If there exist such curves lying on the oriented surfaces $S$ and $S_1$, respectively, we call the surface pair $\{S, S_1\}$ as Bertrand pair surfaces (Fig. 1).

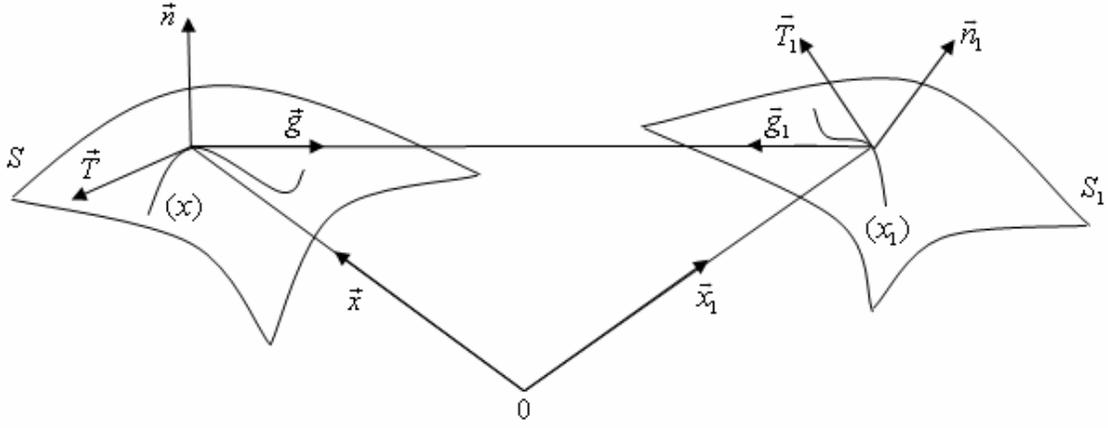

**Fig. 1** Bertrand partner $D$-curves

By considering the Lorentzian casual characters of the surfaces and the curves, from Definition 4.1, it is easily seen that there are five different types of the Bertrand $D$-curves in Minkowski 3-space. Let the pair $\{x, x_1\}$ be a Bertrand $D$-pair. Then according to the character of the surface $S$ we have the followings:

If both the surface $S$ and the curve $x(s)$ lying on $S$ are spacelike then, there are two cases; first one is that both the surface $S_1$ and the curve $x_1(s_1)$ fully lying on $S_1$ are spacelike. In this case we say that the pair $\{x, x_1\}$ is a Bertrand $D$-pair of the type 1. The second case is that both the surface $S_1$ and the curve $x_1(s_1)$ fully lying on $S_1$ are timelike. Then the pair $\{x, x_1\}$ is called a Bertrand $D$-pair of the type 2. If both the surface $S$ and the curve $x(s)$ lying on $S$ are timelike then, there are two cases; one is that both the surface $S_1$ and the curve $x_1(s_1)$ fully lying on $S_1$ are timelike. In this case we say that the pair $\{x, x_1\}$ is a Bertrand $D$-pair of the type 3. The other case is that both the surface $S_1$ and the curve $x_1(s_1)$ fully lying on $S_1$ are spacelike then the pair $\{x, x_1\}$ is a Bertrand $D$-pair of the type 4. If the surface $S$ is timelike and the curve $x(s)$ lying on $S$ is spacelike then the surface $S_1$ is timelike and the curve $x_1(s_1)$ fully lying on $S_1$ is spacelike. In this case we say that the pair $\{x, x_1\}$ is a Bertrand $D$-pair of the type 5.

**Theorem 4.1.** Let $S$ be an oriented surface and $x(s)$ be a Bertrand $D$-curve in $E_1^3$ with arc length parameter $s$ fully lying on $S$. If $S_1$ is another oriented surface and $x_1(s_1)$ is a curve with arc length parameter $s_1$ fully lying on $S_1$, then $x_1(s_1)$ is Bertrand partner $D$-curve of $x(s)$ if and only if the normal curvature $k_n$ of $x(s)$ and the geodesic curvature $k_{g_1}$, the normal curvature $k_{n_1}$ and the geodesic torsion $\tau_{g_1}$ of $x_1(s_1)$ satisfy the following equation

  i) if the pair $\{x, x_1\}$ is a Bertrand $D$-pair of the type 1, then

$$\dot{\tau}_{g_1} = \frac{1}{\lambda}\left[\left(\frac{(1-\lambda k_{g_1})^2 - \lambda^2 \tau_{g_1}^2}{(1-\lambda k_{g_1})}\right)\left(-k_{n_1} + k_n \frac{1-\lambda k_{g_1}}{\cosh\theta}\right) - \frac{\lambda^2 \tau_{g_1} \dot{k}_{g_1}}{1-\lambda k_{g_1}}\right]$$

  ii) if the pair $\{x, x_1\}$ is a Bertrand $D$-pair of the type 2, then

$$\dot{\tau}_{g_1} = \frac{-1}{\lambda}\left[\left(\frac{(1+\lambda k_{g_1})^2 - \lambda^2 \tau_{g_1}^2}{(1+\lambda k_{g_1})}\right)\left(-k_{n_1} + k_n \frac{1+\lambda k_{g_1}}{\sinh\theta}\right) + \left(\frac{-\lambda^2 \dot{k}_{g_1} \tau_{g_1}}{1+\lambda k_{g_1}}\right)\right]$$

***iii)*** *if the pair* $\{x, x_1\}$ *is a Bertrand D-pair of the type 3, then*

$$\dot{\tau}_{g_1} = \frac{-1}{\lambda}\left[\left(\frac{(1+\lambda k_{g_1})^2 - \lambda^2 \tau_{g_1}^2}{(1+\lambda k_{g_1})}\right)\left(-k_{n_1} + k_n \frac{1+\lambda k_{g_1}}{\cosh \theta}\right) - \frac{\lambda^2 \tau_{g_1} \dot{k}_{g_1}}{1+\lambda k_{g_1}}\right]$$

***iv)*** *if the pair* $\{x, x_1\}$ *is a Bertrand D-pair of the type 4, then*

$$\dot{\tau}_{g_1} = \frac{1}{\lambda}\left[\left(\frac{(1-\lambda k_{g_1})^2 - \lambda^2 \tau_{g_1}^2}{(1-\lambda k_{g_1})}\right)\left(-k_{n_1} + k_n \frac{1-\lambda k_{g_1}}{\sinh \theta}\right) - \frac{\lambda^2 \tau_{g_1} \dot{k}_{g_1}}{1-\lambda k_{g_1}}\right]$$

***v)*** *if the pair* $\{x, x_1\}$ *is a Bertrand D-pair of the type 5, then*

$$\dot{\tau}_{g_1} = \frac{1}{\lambda}\left[\left(\frac{(1+\lambda k_{g_1})^2 + \lambda^2 \tau_{g_1}^2}{(1+\lambda k_{g_1})}\right)\left(k_{n_1} - k_n \frac{1+\lambda k_{g_1}}{\cos \theta}\right) + \frac{\lambda^2 \tau_{g_1} \dot{k}_{g_1}}{1+\lambda k_{g_1}}\right]$$

*for some nonzero constants* $\lambda$, *where* $\theta$ *is the angle between the tangent vectors* $T$ *and* $T_1$ *at the corresponding points of* $x$ *and* $x_1$.

**Proof: i)** Suppose that the pair $\{x, x_1\}$ is a Bertrand $D$-pair of the type 1. Denote the Darboux frames of $x(s)$ and $x_1(s_1)$ by $\{\vec{T}, \vec{g}, \vec{n}\}$ and $\{\vec{T}_1, \vec{g}_1, \vec{n}_1\}$, respectively. Then by the definition we can assume that

$$x(s_1) = x_1(s_1) + \lambda(s_1)\vec{g}_1(s_1), \tag{4}$$

for some function $\lambda(s_1)$. By taking derivative of (4) with respect to $s_1$ and applying the Darboux formulas (1) we have

$$\vec{T}\frac{ds}{ds_1} = (1 - \lambda k_{g_1})\vec{T}_1 + \dot{\lambda}\vec{g}_1 + \lambda \tau_{g_1}\vec{n}_1 \tag{5}$$

Since the direction of $\vec{g}_1$ coincides with the direction of $\vec{g}$, i.e., the tangent vector $\vec{T}$ of the curve lies on the plane spanned by the vectors $\vec{T}_1$ and $\vec{n}_1$, we get

$$\dot{\lambda}(s_1) = 0.$$

This means that $\lambda$ is a nonzero constant. Thus, the equality (5) can be written as follows

$$\vec{T}\frac{ds}{ds_1} = (1 - \lambda k_{g_1})\vec{T}_1 + \lambda \tau_{g_1}\vec{n}_1. \tag{6}$$

Furthermore, we have

$$\vec{T} = \cosh \theta \vec{T}_1 + \sinh \theta \vec{n}_1, \tag{7}$$

where $\theta$ is the angle between the tangent vectors $\vec{T}$ and $\vec{T}_1$ at the corresponding points of $x$ and $x_1$. By differentiating this last equation with respect to $s_1$, we get

$$(k_g \vec{g} + k_n \vec{n})\frac{ds}{ds_1} = (\dot{\theta} + k_{n_1})\sinh \theta \vec{T}_1 + (k_{g_1}\cosh \theta + \tau_{g_1}\sinh \theta)\vec{g}_1 + (\dot{\theta} + k_{n_1})\cosh \theta \vec{n}_1. \tag{8}$$

From this equation and the fact that

$$\vec{n} = \sinh \theta \vec{T}_1 + \cosh \theta \vec{n}_1, \tag{9}$$

we get

$$(k_n \sinh \theta \vec{T}_1 + k_g \vec{g} + k_n \cosh \theta \vec{n}_1)\frac{ds}{ds_1} = (\dot{\theta} + k_{n_1})\sinh \theta \vec{T}_1 + (k_{g_1}\cosh \theta + \tau_{g_1}\sinh \theta)\vec{g}_1$$
$$+ (\dot{\theta} + k_{n_1})\cosh \theta \vec{n}_1 \tag{10}$$

Since the direction of $\vec{g}_1$ is coincident with $\vec{g}$ we have

$$\dot{\theta} = -k_{n_1} + k_n \frac{ds}{ds_1}. \tag{11}$$

From (6) and (7) and notice that $\vec{T}_1$ is orthogonal to $\vec{g}_1$ we obtain

$$\frac{ds}{ds_1} = \frac{1-\lambda k_{g_1}}{\cosh\theta} = \frac{\lambda \tau_{g_1}}{\sinh\theta}. \tag{12}$$

Equality (12) gives us

$$\tanh\theta = \frac{\lambda \tau_{g_1}}{1-\lambda k_{g_1}}. \tag{13}$$

By taking the derivative of this equation and applying (11) we get

$$\dot{\tau}_{g_1} = \frac{-1}{\lambda}\left[\left(\frac{(1-\lambda k_{g_1})^2 - \lambda^2 \tau_{g_1}^2}{(1-\lambda k_{g_1})}\right)\left(-k_{n_1} + k_n\frac{1-\lambda k_{g_1}}{\cosh\theta}\right) - \frac{\lambda^2 \tau_{g_1} \dot{k}_{g_1}}{1-\lambda k_{g_1}}\right], \tag{14}$$

that is desired.

Conversely, assume that the equation (14) holds for some nonzero constants $\lambda$. Then by using (12) and (13), (14) gives us

$$k_n\left(\frac{ds}{ds_1}\right)^3 = \lambda\dot{\tau}_{g_1}(1-\lambda k_{g_1}) + \lambda^2 \tau_{g_1} \dot{k}_{g_1} + \left((1-\lambda k_{g_1})^2 - \lambda^2 \tau_{g_1}^2\right)k_{n_1} \tag{15}$$

Let define a curve

$$x(s_1) = x_1(s_1) + \lambda(s_1)\vec{g}_1(s_1). \tag{16}$$

We will prove that $x$ is a Bertrand $D$-curve and $x_1$ is the Bertrand partner $D$-curve of $x$. By taking the derivative of (16) with respect to $s_1$ twice, we get

$$\vec{T}\frac{ds}{ds_1} = (1-\lambda k_{g_1})\vec{T}_1 + \lambda \tau_{g_1}\vec{n}_1, \tag{17}$$

and

$$(k_g \vec{g} + k_n \vec{n})\left(\frac{ds}{ds_1}\right)^2 + \vec{T}\frac{d^2s}{ds_1^2} = (-\lambda\dot{k}_{g_1} + \lambda \tau_{g_1} k_{n_1})\vec{T}_1 + \left((1-\lambda k_{g_1})k_{g_1} + \lambda \tau_{g_1}^2\right)\vec{g}_1 \tag{18}$$
$$+ \left((1-\lambda k_{g_1})k_{n_1} + \lambda\dot{\tau}_{g_1}\right)\vec{n}_1$$

respectively. Taking the cross product of (17) with (18) we have

$$\left[k_g \vec{n} + k_n \vec{g}\right]\left(\frac{ds}{ds_1}\right)^3 = \left[\lambda \tau_{g_1} k_{g_1}(1-\lambda k_{g_1}) + \lambda^2 \tau_{g_1}^3\right]\vec{T}_1$$
$$+ \left[(1-\lambda k_{g_1})^2 k_{n_1} + \lambda\dot{\tau}_{g_1}(1-\lambda k_{g_1}) + \lambda^2 \tau_{g_1}\dot{k}_{g_1} - \lambda^2 \tau_{g_1}^2 k_{n_1}\right]\vec{g}_1 \tag{19}$$
$$+ \left[k_{g_1}(1-\lambda k_{g_1})^2 + \lambda \tau_{g_1}^2(1-\lambda k_{g_1})\right]\vec{n}_1$$

By substituting (15) in (19) we get

$$\left[k_g \vec{n} + k_n \vec{g}\right]\left(\frac{ds}{ds_1}\right)^3 = \left(\lambda \tau_{g_1} k_{g_1}(1-\lambda k_{g_1}) + \lambda^2 \tau_{g_1}^3\right)\vec{T}_1 + k_n\left(\frac{ds}{ds_1}\right)^3\vec{g}_1 \tag{20}$$
$$+ \left(k_{g_1}(1-\lambda k_{g_1})^2 + \lambda \tau_{g_1}^2(1-\lambda k_{g_1})\right)\vec{n}_1$$

Taking the cross product of (17) with (20) we have

$$\left[k_g \vec{g} + k_n \vec{n}\right]\left(\frac{ds}{ds_1}\right)^4 = k_n\left(\frac{ds}{ds_1}\right)^3 \lambda \tau_{g_1}\vec{T}_1 + \left((1-\lambda k_{g_1})^2 - \lambda^2 \tau_{g_1}^2\right)\left(\lambda \tau_{g_1}^2 + k_{g_1}(1-\lambda k_{g_1})\right)\vec{g}_1 \tag{21}$$
$$+ k_n\left(\frac{ds}{ds_1}\right)^3(1-\lambda k_{g_1})\vec{n}_1$$

From (20) and (21) we have

$$(k_g^2 - k_n^2)\left(\frac{ds}{ds_1}\right)^4 n = \left[\lambda k_g k_{g_1} \tau_{g_1}(1-\lambda k_{g_1})\frac{ds}{ds_1} + \lambda^2 k_g \tau_{g_1}^3 \frac{ds}{ds_1} - \lambda \tau_{g_1} k_n^2 \left(\frac{ds}{ds_1}\right)^3\right] T_1$$

$$+ k_n \left(\frac{ds}{ds_1}\right)^2 \left[k_g \left(\frac{ds}{ds_1}\right)^2 - \lambda \tau_{g_1}^2 - k_{g_1}(1-\lambda k_{g_1})\right] g_1 \qquad (22)$$

$$+ \left[k_g k_{g_1}(1-\lambda k_{g_1})^2 \frac{ds}{ds_1} + \lambda \tau_{g_1}^2 k_g (1-\lambda k_{g_1})\frac{ds}{ds_1} - k_n^2 (1-\lambda k_{g_1})\left(\frac{ds}{ds_1}\right)^3\right] n_1$$

Furthermore, from (17) and (20) we get

$$\begin{cases} \left(\dfrac{ds}{ds_1}\right)^2 = (1-\lambda k_{g_1})^2 - \lambda^2 \tau_{g_1}^2, \\ k_g \left(\dfrac{ds}{ds_1}\right)^2 = k_{g_1}(1-\lambda k_{g_1}) + \lambda \tau_{g_1}^2 \end{cases} \qquad (23)$$

respectively. Substituting (23) in (22) we obtain

$$(k_g^2 - k_n^2)\left(\frac{ds}{ds_1}\right)^4 n = \left[\lambda k_g k_{g_1} \tau_{g_1}(1-\lambda k_{g_1})\frac{ds}{ds_1} + \lambda^2 k_g \tau_{g_1}^3 \frac{ds}{ds_1} - \lambda \tau_{g_1} k_n^2 \left(\frac{ds}{ds_1}\right)^3\right] T_1$$

$$+ \left[k_g k_{g_1}(1-\lambda k_{g_1})^2 \frac{ds}{ds_1} + \lambda \tau_{g_1}^2 k_g (1-\lambda k_{g_1})\frac{ds}{ds_1} - k_n^2 (1-\lambda k_{g_1})\left(\frac{ds}{ds_1}\right)^3\right] n_1 \qquad (24)$$

Equality (17) and (24) shows that the vectors $\vec{T}$ and $\vec{n}$ lie on the plane $sp\{\vec{T}_1, \vec{n}_1\}$. So, at the corresponding points of the curves, the Darboux frame element $\vec{g}$ of $x$ coincides with the Darboux frame element $\vec{g}_1$ of $x_1$, i.e, the curves $x$ and $x_1$ are Bertrand $D$-pair curves.

Let now give the characterizations of Bertrand partner $D$-curves of the type 1 in some special cases. Assume that $x(s)$ be an asymptotic line. Then, from (14) we have the following special cases:

**i)** Consider that $x_1(s_1)$ is a geodesic curve. Then $x_1(s_1)$ is Bertrand partner $D$-curve of $x(s)$ if and only if the following equation holds,

$$\lambda \dot{\tau}_{g_1} = k_{n_1}\left(1 - \lambda^2 \tau_{g_1}^2\right)$$

**ii)** Assume that $x_1(s_1)$ is also an asymptotic line. Then $x_1(s_1)$ is Bertrand partner $D$-curve of $x(s)$ if and only if the geodesic torsion $\tau_{g_1}$ of $x_1(s_1)$ satisfies the following equation,

$$\dot{\tau}_{g_1} = \frac{\lambda \tau_{g_1} \dot{k}_{g_1}}{1 - \lambda k_{g_1}}.$$

In this case, the Frenet frame of the curve $x_1(s_1)$ coincides with its Darboux frame. From (2) we have $k_{g_1} = \kappa_1$ and $\tau_{g_1} = \tau_1$. So, the Bertrand partner $D$-curves become the Bertrand partner curves, i.e., if both $x(s)$ and $x_1(s_1)$ are asymptotic lines then, the definition and the characterizations of the Bertrand partner $D$-curves involve those of the Bertrand partner curves in Minkowski 3-space. Then, a new characterization of Bertrand curves can be given as follows

***Corollary 4.1.*** *Let $x(s)$ be a spacelike Bertrand curve with arclength parameter $s$ in Minkowski 3-space $E_1^3$. Then the spacelike curve $x_1(s_1)$ is a Bertrand partner curve of $x(s)$ if and only if the curvature $\kappa_1$ and the torsion $\tau_1$ of $x_1(s_1)$ satisfy the following equation*

$$\dot{\tau}_{g_1} = \frac{\lambda \tau_{g_1} \dot{k}_{g_1}}{1 - \lambda k_{g_1}}$$

*for some nonzero constants $\lambda$.*

iii) If $x_1(s_1)$ is a principal line then $x_1(s_1)$ is Bertrand partner $D$-curve of $x(s)$ if and only if the geodesic curvature $k_{g_1}$ and the geodesic torsion $\tau_{g_1}$ of $x_1(s_1)$ satisfy the following equality,

$$\lambda \dot{\tau}_{g_1} = k_{n_1}(1 - \lambda k_{g_1}).$$

The proofs of the statement (ii), (iii), (iv) and (v) of Theorem 4.1 and the particular cases given above can be given by the same way of the proof of statement (i).

***Theorem 4.2.*** *Let the pair $\{x, x_1\}$ is a Bertrand $D$-pair of the type 1. Then the relation between geodesic curvature $k_g$, geodesic torsion $\tau_g$ of $x(s)$ and the geodesic curvature $k_{g_1}$, the geodesic torsion $\tau_{g_1}$ of $x_1(s_1)$ is given as follows*

  i) *if the pair $\{x, x_1\}$ is a Bertrand $D$-pair of the type 1, then*
  $$k_g - k_{g_1} = \lambda(k_g k_{g_1} + \tau_g \tau_{g_1})$$
  ii) *if the pair $\{x, x_1\}$ is a Bertrand $D$-pair of the type 2, then*
  $$k_g + k_{g_1} = -\lambda(k_g k_{g_1} + \tau_g \tau_{g_1})$$
  iii) *if the pair $\{x, x_1\}$ is a Bertrand $D$-pair of the type 3, then*
  $$k_g - k_{g_1} = -\lambda(k_g k_{g_1} - \tau_g \tau_{g_1})$$
  iv) *if the pair $\{x, x_1\}$ is a Bertrand $D$-pair of the type 4, then*
  $$k_g + k_{g_1} = \lambda(k_g k_{g_1} + \tau_g \tau_{g_1})$$
  v) *if the pair $\{x, x_1\}$ is a Bertrand $D$-pair of the type 4, then*
  $$k_g - k_{g_1} = \lambda(\tau_g \tau_{g_1} - k_g k_{g_1}).$$

**Proof: i)** Suppose that the pair $\{x, x_1\}$ is a Bertrand $D$-pair of the type 1. Then by definition from (16) we can write
$$x_1(s_1) = x(s_1) - \lambda(s_1)\vec{g}(s_1) \tag{23}$$
for some constants $\lambda$. By differentiating (23) with respect to $s_1$ we have
$$\vec{T}_1 = (1 + \lambda k_g)\frac{ds}{ds_1}\vec{T} - \lambda \tau_g \frac{ds}{ds_1}\vec{n} \tag{24}$$
By the definition we have
$$\vec{T}_1 = \cosh\theta \vec{T} - \sinh\theta \vec{n} \tag{25}$$
From (24) and (25) we obtain
$$\cosh\theta = (1 + \lambda k_g)\frac{ds}{ds_1}, \quad \sinh\theta = \lambda \tau_g \frac{ds}{ds_1} \tag{26}$$
Using (12) and (26) it is easily seen that
$$k_g - k_{g_1} = \lambda(k_g k_{g_1} + \tau_g \tau_{g_1}).$$

From Theorem 4.2, we obtain the following special cases.
Let the pair $\{x, x_1\}$ be a Bertrand $D$-pair of the type 1. Then,

    **i)** if one of the curves $x$ and $x_1$ is a principal line, then the relation between the geodesic curvatures $k_g$ and $k_{g_1}$ is

$$k_g - k_{g_1} = \lambda k_g k_{g_1}$$

    **ii)** if $x_1$ is a geodesic curve, then the geodesic curvature of the curve $x$ is given by

$$k_g = \lambda \tau_g \tau_{g_1}$$

    **iii)** if $x$ is a geodesic curve, then the geodesic curvature of the curve $x_1$ is given by

$$k_{g_1} = -\lambda \tau_g \tau_{g_1}$$

**Theorem 4.3.** *Let $\{x, x_1\}$ be Bertrand $D$-pair of the type 1. Then the following relations hold:*

    **i)** $k_{n_1} = k_n \dfrac{ds}{ds_1} - \dfrac{d\theta}{ds_1}$

    **ii)** $\tau_g \dfrac{ds}{ds_1} = k_{g_1} \sinh\theta - \tau_{g_1} \cosh\theta$

    **iii)** $k_g \dfrac{ds}{ds_1} = k_{g_1} \cosh\theta + \tau_{g_1} \sinh\theta$

    **iv)** $\tau_{g_1} = (-k_g \sinh\theta + \tau_g \cosh\theta)\dfrac{ds}{ds_1}$

**Proof: i)** Since the pair $\{x, x_1\}$ is a Bertrand $D$-pair of the type 1, we have $\langle \vec{T}, \vec{T_1} \rangle = \cosh\theta$.
By differentiating this equation with respect to $s_1$ we have

$$\left\langle (k_g \vec{g} + k_n \vec{n})\dfrac{ds}{ds_1}, \vec{T_1} \right\rangle + \left\langle \vec{T}, k_{g_1} \vec{g_1} + k_{n_1} \vec{n_1} \right\rangle = \sinh\theta \dfrac{d\theta}{ds_1}.$$

Using the fact that the direction of $g_1$ coincides with the direction of $g$ and

$$\vec{T_1} = \cosh\theta \vec{T} - \sinh\theta \vec{n}, \quad \vec{n_1} = -\sinh\theta \vec{T} + \cosh\theta \vec{n} \qquad (27)$$

we easily get that

$$k_{n_1} = k_n \dfrac{ds}{ds_1} - \dfrac{d\theta}{ds_1}.$$

**ii)** By definition we have $\langle \vec{n}, \vec{g_1} \rangle = 0$. Differentiating this equation with respect to $s_1$ we have

$$\left\langle (k_n \vec{T} + \tau_g \vec{g})\dfrac{ds}{ds_1}, \vec{g_1} \right\rangle + \left\langle \vec{n}, -k_{g_1} \vec{T_1} + \tau_{g_1} \vec{n_1} \right\rangle = 0.$$

By (27) we obtain

$$\tau_g \dfrac{ds}{ds_1} = k_{g_1} \sinh\theta - \tau_{g_1} \cosh\theta$$

**iii)** By differentiating the equation $\langle \vec{T}, \vec{g_1} \rangle = 0$ with respect to $s_1$ we get

$$\left\langle (k_g \vec{g} + k_n \vec{n})\dfrac{ds}{ds_1}, \vec{g_1} \right\rangle + \left\langle \vec{T}, (-k_{g_1} \vec{T_1} + \tau_{g_1} \vec{n_1}) \right\rangle = 0.$$

From (27) it follows that

$$k_g \dfrac{ds}{ds_1} = k_{g_1} \cosh\theta + \tau_{g_1} \sinh\theta.$$

**iv)** By differentiating the equation $\langle \vec{n}_1, \vec{g} \rangle = 0$ with respect to $s_1$ we obtain

$$\langle (k_{n_1}\vec{T}_1 + \tau_{g_1}\vec{g}_1), \vec{g} \rangle + \langle \vec{n}_1, (-k_g\vec{T} + \tau_g\vec{n}) \frac{ds}{ds_1} \rangle = 0,$$

and using the fact that direction of $g_1$ coincides with the direction of $g$ and

$$\vec{T} = \cosh\theta \vec{T}_1 + \sinh\theta \vec{n}_1, \quad \vec{n} = \sinh\theta \vec{T}_1 + \cosh\theta \vec{n}_1$$

we get

$$\tau_{g_1} = (-k_g \sinh\theta + \tau_g \cosh\theta) \frac{ds}{ds_1}.$$

The statements of Theorem 4.3 for the pairs $\{x, x_1\}$ of the type 2, 3, 4, and 5 can be given as follows and the proofs can be easily done by the same way of the case the pairs $\{x, x_1\}$ is of the type 1.

| For the pair $\{x, x_1\}$ of the type 2 | For the pair $\{x, x_1\}$ of the type 3 |
|---|---|
| i) $k_{n_1} = \dfrac{d\theta}{ds_1} + k_n \dfrac{ds}{ds_1}$ | i) $k_{n_1} = k_n \dfrac{ds}{ds_1} + \dfrac{d\theta}{ds_1}$ |
| ii) $\tau_g \dfrac{ds}{ds_1} = k_{g_1} \cosh\theta - \tau_{g_1} \sinh\theta$ | ii) $\tau_g \dfrac{ds}{ds_1} = k_{g_1} \sinh\theta - \tau_{g_1} \cosh\theta$ |
| iii) $k_g \dfrac{ds}{ds_1} = k_{g_1} \sinh\theta + \tau_{g_1} \cosh\theta$ | iii) $k_g \dfrac{ds}{ds_1} = k_{g_1} \cosh\theta + \tau_{g_1} \sinh\theta$ |
| iv) $\tau_{g_1} = (k_g \cosh\theta - \tau_g \sinh\theta) \dfrac{ds}{ds_1}$ | iv) $\tau_{g_1} = (-k_g \sinh\theta + \tau_g \cosh\theta) \dfrac{ds}{ds_1}$ |
| For the pair $\{x, x_1\}$ of the type 4 | For the pair $\{x, x_1\}$ of the type 5 |
| i) $k_{n_1} = k_n \dfrac{ds}{ds_1} - \dfrac{d\theta}{ds_1}$ | i) $k_{n_1} = k_n \dfrac{ds}{ds_1} + \dfrac{d\theta}{ds_1}$ |
| ii) $\tau_g \dfrac{ds}{ds_1} = k_{g_1} \cosh\theta + \tau_{g_1} \sinh\theta$ | ii) $\tau_g \dfrac{ds}{ds_1} = -k_{g_1} \sin\theta + \tau_{g_1} \cos\theta$ |
| iii) $k_g \dfrac{ds}{ds_1} = k_{g_1} \sinh\theta + \tau_{g_1} \cosh\theta$ | iii) $k_g \dfrac{ds}{ds_1} = k_{g_1} \cos\theta + \tau_{g_1} \sin\theta$ |
| iv) $\tau_{g_1} = (k_g \cosh\theta - \tau_g \sinh\theta) \dfrac{ds}{ds_1}$ | iv) $\tau_{g_1} = (k_g \sin\theta + \tau_g \cos\theta) \dfrac{ds}{ds_1}$ |

Let now $x$ be a Bertrand $D$-curve and $x_1$ be a Bertrand partner $D$-curve of $x$ and the pair $\{x, x_1\}$ be of the type 1. From the first equation of (3) and by using the fact that $\vec{n}_1 = -\sinh\theta \vec{T} + \cosh\theta \vec{n}$ we have

$$k_{g_1} = \left[(1 + \lambda k_g)\cosh\theta - \lambda \tau_g \sinh\theta\right]\left[-k_g - \lambda k_g^2 + \lambda \tau_g^2\right]\left(\frac{ds}{ds_1}\right)^3. \tag{28}$$

Then we can give the following corollary.

**Corollary 4.2.** *Let $x$ be a Bertrand $D$-curve and $x_1$ be a Bertrand partner $D$-curve of $x$ and the pair $\{x, x_1\}$ be of the type 1. Then the relations between the geodesic curvature $k_{g_1}$ of $x_1(s_1)$ and the geodesic curvature $k_g$ and the geodesic torsion $\tau_g$ of $x(s)$ are given as follows.*

*i) If $x$ is a geodesic curve, then the geodesic curvature $k_{g_1}$ of $x_1(s_1)$ is*

$$k_{g_1} = \lambda \tau_g^2 \left(\frac{ds}{ds_1}\right)^3 (\cosh\theta - \lambda\tau_g \sinh\theta). \tag{29}$$

*ii) If $x$ is a principal line, then the relation between the geodesic curvatures $k_{g_1}$ and $k_g$ is given by*

$$k_{g_1} = -(k_g + 2\lambda k_g^2 + \lambda^2 k_g^3)\left(\frac{ds}{ds_1}\right)^3 \cosh\theta. \tag{30}$$

If the pair $\{x, x_1\}$ is of the type 2, 3, 4 or 5 then the geodesic curvature of the curve $x_1(s_1)$ is given as follows

| If the pair $\{x, x_1\}$ is of the type 2 | If the pair $\{x, x_1\}$ is of the type 3 |
|---|---|
| $k_{g_1} = \left[(1+\lambda k_g)\sinh\theta - \lambda\tau_g \cosh\theta\right]$ $\left[k_g + \lambda k_g^2 - \lambda\tau_g^2\right]\left(\frac{ds}{ds_1}\right)^3$ | $k_{g_1} = \left[(1-\lambda k_g)\cosh\theta + \lambda\tau_g \sinh\theta\right]$ $\left[-k_g + \lambda k_g^2 - \lambda\tau_g^2\right]\left(\frac{ds}{ds_1}\right)^3$ |
| If the pair $\{x, x_1\}$ is of the type 4 | If the pair $\{x, x_1\}$ is of the type 5 |
| $k_{g_1} = \left[(1-\lambda k_g)\sinh\theta + \lambda\tau_g \cosh\theta\right]$ $\left[k_g - \lambda k_g^2 + \lambda\tau_g^2\right]\left(\frac{ds}{ds_1}\right)^3$ | $k_{g_1} = \left[(1-\lambda k_g)\cos\theta + \lambda\tau_g \sin\theta\right]$ $\left[-k_g + \lambda k_g^2 + \lambda\tau_g^2\right]\left(\frac{ds}{ds_1}\right)^3$ |

and the statements in Corollary 4.2 are obtained by the same way.

Similarly, From the second equation of (3) and by using the fact that $g$ is coincident with $g_1$, i.e., $\vec{n}_1 = -\sinh\theta \vec{T} + \cosh\theta \vec{n}$, the geodesic torsion $\tau_{g_1}$ of $x_1$ is given by

$$\tau_{g_1} = \left[(\tau_g + \lambda k_g \tau_g)\cosh^2\theta + (-k_g - \lambda k_g^2 + \lambda\tau_g^2)\sinh\theta\cosh\theta + \lambda\tau_g k_g \sinh^2\theta\right]\left(\frac{ds}{ds_1}\right)^2. \tag{31}$$

From (31) we can give the following corollary.

**Corollary 4.3.** *Let $x$ be a Bertrand $D$-curve and $x_1$ be a Bertrand partner $D$-curve of $x$ and the pair $\{x, x_1\}$ be of the type 1. Then the relations between the geodesic torsion $\tau_{g_1}$ of $x_1(s_1)$ and the geodesic curvature $k_g$ and the geodesic torsion $\tau_g$ of $x(s)$ are given as follows.*
*i) If $x$ is a geodesic curve then the geodesic torsion of $x_1$ is*

$$\tau_{g_1} = (\tau_g \cosh^2\theta + \lambda\tau_g^2 \sinh\theta\cosh\theta)\left(\frac{ds}{ds_1}\right)^2. \tag{32}$$

**ii)** If $x$ is a principal line then the relation between $\tau_{g_1}$ and $k_g$ is

$$\tau_{g_1} = -(k_g + \lambda k_g^2)\sinh\theta\cosh\theta\left(\frac{ds}{ds_1}\right)^2. \tag{33}$$

Furthermore, by using (12) and (13), from (32) and (33) we have the following corollary.

**Corollary 4.4. i)** Let $\{x, x_1\}$ be Bertrand $D$-pair of the type 1 and let $x$ be a geodesic line. Then the geodesic torsion $\tau_{g_1}$ of $x_1(s_1)$ is given by

$$\tau_{g_1} = \tau_g(1 - \lambda k_{g_1})\left[(1 - \lambda k_{g_1}) + \lambda^2 \tau_g \tau_{g_1}\right]. \tag{34}$$

**ii)** Let $\{x, x_1\}$ be Bertrand $D$-pair of the type 1 and let $x$ be a principal line. Then the relation between the geodesic curvatures $k_g$ and $k_{g_1}$ is given as follows

$$k_g(1 + \lambda k_g)(1 - \lambda k_{g_1}) = -\frac{1}{\lambda} = \text{constant}. \tag{35}$$

When the pair $\{x, x_1\}$ is of the type 2, 3, 4 or 5, then the relations which give the geodesic torsion $\tau_{g_1}$ of $x_1(s_1)$ are given as follows.

| For the pair $\{x, x_1\}$ of the type 2 | For the pair $\{x, x_1\}$ of the type 3 |
|---|---|
| $\tau_{g_1} = \left[\tau_g \sinh^2\theta - \lambda\tau_g k_g \right.$ $\left. +(\lambda\tau_g^2 - k_g - \lambda k_g^2))\sinh\theta\cosh\theta\right]\left(\frac{ds}{ds_1}\right)^2$ | $\tau_{g_1} = \left[(\tau_g - \lambda k_g \tau_g)\cosh^2\theta \right.$ $+(-k_g + \lambda k_g^2 + \lambda\tau_g^2)\sinh\theta\cosh\theta$ $\left. -\lambda\tau_g k_g \sin^2\theta\right]\left(\frac{ds}{ds_1}\right)^2$ |
| For the pair $\{x, x_1\}$ of the type 4 | For the pair $\{x, x_1\}$ of the type 5 |
| $\tau_{g_1} = \left[(\tau_g + \lambda k_g \tau_g)\sinh^2\theta \right.$ $-(\lambda\tau_g^2 + k_g + \lambda k_g^2)\sinh\theta\cosh\theta$ $\left. +\lambda\tau_g k_g \cosh^2\theta\right]\left(\frac{ds}{ds_1}\right)^2$ | $\tau_{g_1} = \left[(\tau_g - \lambda k_g \tau_g)\cos^2\theta \right.$ $+(k_g - \lambda k_g^2 + \lambda\tau_g^2)\sin\theta\cos\theta$ $\left. +\lambda\tau_g k_g \sin^2\theta\right]\left(\frac{ds}{ds_1}\right)^2$ |

## 4. Conclusions

In this paper, the definition and characterizations of Bertrand partner $D$-curves are given which is a new study of associated curves lying on surfaces. It is shown that the definition and the characterizations of Bertrand partner $D$-curves include those of Bertrand partner curves in some special cases. Furthermore, the relations between the geodesic curvatures, the normal curvatures and the geodesic torsions of these curves are given.